\documentclass{conm-p-l}

\usepackage{graphics}

\newtheorem{theorem}{Theorem}[section]

\theoremstyle{definition}
\newtheorem{definition}[theorem]{Definition}

\theoremstyle{remark}
\newtheorem{remark}[theorem]{Remark}

\numberwithin{equation}{section}

\newcommand{\LL}{{\mathcal L}}
\newcommand{\tla}{\tilde{\lambda}}

\newcommand{\hk}{\hat{k}}
\newcommand{\Z}{\ZZ^2/\{0\}}

\newcommand{\CC}{{\mathbb C}}
\newcommand{\RR}{{\mathbb R}}
\newcommand{\ZZ}{{\mathbb Z}}

\newcommand{\e}{\varepsilon}
\newcommand{\vth}{\vartheta}

\renewcommand{\k}{\kappa}
\newcommand{\ga}{\gamma}
\newcommand{\Ga}{\Gamma}

\newcommand{\dl}{\delta}
\newcommand{\Dl}{\Delta}
\renewcommand{\th}{\theta}
\newcommand{\ra}{\rightarrow}
\newcommand{\al}{\alpha}
\newcommand{\be}{\beta}
\newcommand{\sg}{\sigma}
\newcommand{\Sg}{\Sigma}

\newcommand{\pa}{\partial}
\newcommand{\z}{\zeta}

\newcommand{\La}{\Lambda}

\newcommand{\la}{\lambda}

\newcommand{\nid}{\noindent}

\newcommand{\om}{\omega}
\newcommand{\Om}{\Omega}
\newcommand{\na}{\nabla}

\newcommand{\tq}{\tilde{q}}

\newcommand{\W}{{\mathcal W}}

\def\maprightu#1{\smash{
    \mathop{\longrightarrow}\limits^{#1}}}
\def\maprightd#1{\smash{
    \mathop{\longrightarrow}\limits_{#1}}}
\def\mapdownl#1{
    \llap{$\vcenter{\hbox{$\scriptstyle#1$}}$}\Big\downarrow}
\def\mapdownr#1{\Big\downarrow
    \rlap{$\vcenter{\hbox{$\scriptstyle#1$}}$}}

\begin{document}

\title{Chaos in PDEs and Lax Pairs of Euler Equations}

\author{Yanguang (Charles)  Li}
\address{Department of Mathematics, University of Missouri, 
Columbia, MO 65211}
\curraddr{}
\email{cli@math.missouri.edu}
\thanks{}


\date{}



\begin{abstract}
Recently, the author and collaborators have developed a systematic 
program for proving the existence of homoclinic orbits in partial 
differential equations. Two typical forms of homoclinic orbits thus 
obtained are: (1). transversal homoclinic orbits, (2). Silnikov 
homoclinic orbits. Around the transversal homoclinic orbits in infinite 
dimensional autonomous systems, the author was able to prove the 
existence of chaos through a shadowing lemma. Around the Silnikov 
homoclinic orbits, the author was able to prove the existence of chaos 
through a horseshoe construction. 

Very recently, there has been a breakthrough by the author in finding 
Lax pairs for Euler equations of incompressible inviscid fluids. 
Further results have been obtained by the author and collaborators. 
\end{abstract}

\maketitle








\section{Introduction}

Unlike chaos in finite dimensional systems, the area of chaos in partial 
differential equations has been a long standing open field.
During the past ten years, there have been some remarkable developments 
in the area of chaos in partial differential equations by the author 
and collaborators. Our main results can be briefly described as follows: 
A systematic program is developed for proving the existence of homoclinic 
orbits for perturbed soliton equations. The program involves machineries 
from integrable theory, dynamical systems, and partial differential equations.
The types of homoclinic orbits thus obtained can be either transversal 
or of Silnikov. With regard to transversal homoclinic orbits, shadowing 
lemmas are utilized or developed to prove the existence of chaos. With 
regard to Silnikov homoclinic orbits, Smale horseshoes are constructed 
to prove the existence of chaos. The main machineries from integrable theory 
are the isospectral theory and Darboux transformations. The main machineries 
from dynamical systems are persistence of invariant manifolds, Fenichel 
fibers, Melnikov measurement, and other measurements. The main machineries 
from partial differential equations are local well-posedness and regularity.

Another exciting development happened to Euler equations of incompressible 
inviscid fluids. Lax pairs for both 2D and 3D Euler equations was  
found by the author and Steve Childress. A Darboux transformation for the 
Lax pair of 2D Euler equation was also found by Artyom Yurov and the author. 
The philosophical significance of the existence of Lax pairs for Euler 
equations is fundamentally important. If one defines integrability of an 
equation by the existence of a Lax pair, then both 2D and 3D Euler equations 
are integrable. More importantly, both 2D and 3D Navier-Stokes equations 
at high Reynolds numbers are near integrable systems. Such a point of view 
changes our old ideology on Euler and Navier-Stokes equations.

The most important application of the theory on chaos in partial 
differential equations in theoretical physics will be on the study of 
turbulence. Our approach of studying turbulence is different from many 
other studies in which one starts with Stokes equation to prove results on
Navier-Stokes equations for small Reynolds number. In our studies, we 
start with Euler equations and view Navier-Stokes equations at large 
Reynolds number as singular perturbations of Euler equations. For this goal, 
we chose the 2D Navier-Stokes equations under periodic boundary conditions 
to begin a dynamical system study on 2D turbulence. We started with a 
simple fixed point of 2D Euler equation, and studied the linearized 2D 
Euler equation at the fixed point. A complete spectral theorem is obtained.
In particular, we found unstable eigenvalues. Then, naturally we are 
interested in the corresponding unstable manifold. Such an unstable manifold 
forms the prototype of the attractor of Navier-Stokes equations at large 
Reynolds number. And long term turbulence lives around the attractor. 
Finding the unstable manifold is not successful yet. Partial success is 
accomplished in finding the corresponding unstable manifold for some 
Galerkin truncation. Here the unstable manifold is the 2D surface of 
an ellipsoid. Stable and unstable manifolds together form a lip-shape 
configuration. Whether or not the stable and unstable manifolds for 2D Euler 
equation have the same topology as above is not known. In fact, the 
existence of stable and unstable manifolds for 2D Euler equation has not been 
proved. 

\section{Homoclinic Orbits}

In terms of proving the existence of a homoclinic orbit, 
the most common tool is the so-called Melnikov integral method \cite{Mel63} 
\cite{Arn64}. This method was subsequently developed by Holmes and Marsden 
\cite{GH83}, and most recently by Wiggins \cite{Wig88}. For partial 
differential equations, this method was mainly developed by Li et al.
\cite{LMSW96} \cite{Li01b} \cite{LM97} \cite{Li02f}.

There are two derivations for the Melnikov integrals. One is the 
so-called geometric argument \cite{GH83} \cite{Wig88} \cite{LM97} 
\cite{LMSW96} \cite{Li01b}.
The other is the so-called Liapunov-Schmitt argument \cite{CHMP80} 
\cite{CH82}. The Liapunov-Schmitt argument is a fixed-point type argument 
which directly leads to the existence of a homoclinic orbit. The condition 
for the existence of a fixed point is the Melnikov integral. 
The geometric argument is a signed distance argument which applies 
to more general situations than the Liapunov-Schmitt argument. It turns out 
that the geometric argument is a much more powerful machinary than the 
Liapunov-Schmitt argument. In particular, the geometric argument can handle 
geometric singular perturbation problems. I shall also mention an 
interesting derivation in \cite{Arn64}.

In establishing the existence of homoclinic orbits in high dimensions, one 
often needs other tools besides the Melnikov analysis. For example, when 
studying orbits 
homoclinic to fixed points created through resonances in 
($n \geq 4$)-dimensional near-integrable systems, one often needs tools 
like Fenichel fibers, as presented in previous chapter, to set up geometric 
measurements for locating such 
homoclinic orbits. Such homoclinic orbits often have a geometric singular 
perturbation nature. In such cases, the Liapunov-Schmitt argument can not 
be applied. For such works on finite dimensional systems, see for 
example \cite{Kov92a} \cite{Kov92b} \cite{LM97}.
For such works on infinite dimensional systems, see for example 
\cite{LMSW96} \cite{Li01b}.

\subsection{Silnikov Homoclinic Orbits in Nonlinear 
Schr\"odinger (NLS) Equation  Under Regular 
Perturbations \label{horrnls}}

Consider the regularly perturbed nonlinear Schr\"odinger (NLS) equation 
\cite{LMSW96},
\begin{equation}
iq_t = q_{xx} +2 [ |q|^2 - \om^2] q +i \e [\hat{\pa}^2_xq - \al q +\be ] \ ,
\label{rpnls}
\end{equation}
where $q = q(t,x)$ is a complex-valued function of the two real 
variables $t$ and $x$, $t$ represents time, and $x$ represents
space. $q(t,x)$ is subject to periodic boundary condition of period 
$2 \pi$, and even constraint, i.e., 
\[
q(t,x + 2 \pi) = q(t,x)\ , \ \ q(t,-x) = q(t,x)\ .
\]
$\om$ is a positive constant, $\al >0$ and $\be >0$ are constants, 
$\hat{\pa}^2_x$ is a bounded Fourier multiplier,
\[
\hat{\pa}^2_x q = -\sum_{k=1}^{N}k^2 \xi_k \tq_k \cos kx\ ,
\]
$\xi_k = 1$ when $k \leq N$, $\xi_k = 8k^{-2}$ when $k>N$, for some 
fixed large $N$, and $\e > 0$ is the perturbation parameter.
The following theorem was proved in \cite{LMSW96}.
\begin{theorem}
There exists a $\e_0 > 0$, such that for any $\e \in (0, \e_0)$, there 
exists a codimension 1 surface in the external parameter space 
$(\alpha,\beta, \om) \in \RR^+\times \RR^+\times \RR^+$ where $\om \in 
(\frac{1}{2}, 1)$, and $\al \om < \be$. For any 
$(\alpha ,\beta, \omega)$ on the codimension 1
surface, the regularly perturbed nonlinear Schr\"odinger equation 
(\ref{rpnls}) possesses a symmetric pair of Silnikov homoclinic orbits 
asymptotic 
to a saddle $Q_\epsilon$. The codimension 1 surface has the approximate 
representation given by $\al = 1/\k(\om)$, where $\k(\om)$ is plotted 
in Figure \ref{kappa}.
\label{rhorbit}
\end{theorem}
\begin{figure}
\caption{The graph of $\k(\om)$.}
\label{kappa}
\end{figure}
To prove the theorem, one starts from the invariant plane 
\[
\Pi=\{ q\mid \ \partial_x q=0 \}.
\]
On $\Pi$, there is a saddle $Q_\e =\sqrt{I} e^{i\th}$ to which 
the symmetric pair of Silnikov homoclinic orbits will be asymptotic to, where
\begin{equation}
I=\omega^2-\epsilon \frac{1}{2\omega}\sqrt{\beta^2-\alpha^2\omega^2}+\cdots ,
\quad \cos \theta  =\frac{\alpha \sqrt{I}}{\beta}, \quad \theta \in 
(0,\frac{\pi}{2}).
\label{Qec}
\end{equation}
Its eigenvalues are
\begin{equation}
\la_n^\pm = -\e [\al +\xi_n n^2]\pm 2 \sqrt{(\frac{n^2}{2} +
\om^2-I)(3I -\om^2 -\frac{n^2}{2} )}\ , 
\label{Qev}
\end{equation}
where $n=0,1,2, \cdots $, $\om \in (\frac{1}{2}, 1)$, $\xi_n = 1$ when 
$n \leq N$, $\xi_n = 8n^{-2}$ when $n>N$, for some fixed large $N$,
and $I$ is given in 
(\ref{Qec}). The crucial points to notice are: (1). only $\la_0^+$ and 
$\la_1^+$ have positive real parts, $\mbox{Re}\{ \la_0^+\} < \mbox{Re}\{ 
\la_1^+\} $; (2). all the other eigenvalues have negative real parts among 
which the absolute value of $\mbox{Re}\{ \la_2^+\}=\mbox{Re}\{ \la_2^-\}$ 
is the smallest; (3). $|\mbox{Re}\{ \la_2^+\}| < \mbox{Re}\{ \la_0^+\}$. 
Actually, items (2) and (3) are the main characteristics of Silnikov 
homoclinic orbits. 

The unstable manifold $W^u(Q_\e)$ of $Q_\e$ has a fiber representation 
\cite{LMSW96}. The Melnikov measurement measures the 
signed distance between $W^u(Q_\e)$ and a center-stable manifold 
$W^{cs}_\e$ \cite{LMSW96}. By virtue of a Fenichel Fiber Theorem 
\cite{LMSW96}, one can show that, to the leading order in
$\e$, the signed distance is given by the Melnikov integral
\begin{eqnarray*}
M &=& \int^{+\infty}_{-\infty}\int^{2\pi}_0 
 [\partial_qF_1(q_0(t))(\hat{\partial}^2_xq_0(t)-\alpha
q_0(t)+\beta ) \\
& & \quad \quad + \partial_{\bar{q}}F_1(q_0(t))
(\hat{\partial}^2_x\overline{q_0(t)}-\alpha 
\overline{q_0(t)}+\beta )]dxdt,
\end{eqnarray*}
where $q_0(t)$ is a homoclinic orbit for the integrable NLS, which 
can be generated through Darboux transformations \cite{LM94};
and $\partial_qF_1$ and $\partial_{\bar{q}}F_1$ are Melnikov vectors 
which can be generated through Floquet discriminant in the isospectral 
theory of the integrable NLS \cite{LM94}. The zero of the 
signed distance implies the existence of an orbit in $W^u(Q_\e)\cap 
W^{cs}_\e$. The stable manifold $W^s(Q_\e)$ of $Q_\e$ is a codimension 
1 submanifold in $W^{cs}_\e$. To locate a homoclinic orbit, one needs 
to set up a second measurement measuring the signed distance between 
the orbit in $W^u(Q_\e)\cap W^{cs}_\e$ and $W^s(Q_\e)$ inside $W^{cs}_\e$.
To set up this signed distance, first one can rather easily track 
the (perturbed) orbit by an unperturbed orbit to an $O (\e)$ neighborhood 
of $\Pi$, then one needs to prove the size of $W^s(Q_\e)$ to be $O (\e^\nu)$
($\nu <1$) with normal form transform. To the leading order in
$\e$, the zero of the second signed distance is given by
\[
\be \cos \ga = \frac {\al \om (\Dl \ga )} {2 \sin \frac {\Dl \ga }{2}} \ ,
\]
where $\Dl \ga = -4 \vth_0$ and $\vth_0$ is a phase shift of $q_0(t)$.
To the leading order in $\e$, the common zero of the two second signed 
distances satisfies $\al = 1/\k(\om)$, where $\k(\om)$ is plotted 
in Figure \ref{kappa}. Then the claim of the theorem is proved by virtue 
of the implicit function theorem. For rigorous details, see \cite{LMSW96}.

\subsection{Silnikov Homoclinic Orbits in NLS Under Singular 
Perturbations \label{horsnls}}

Consider the singularly perturbed nonlinear Schr\"odinger equation 
\cite{Li01b},
\begin{equation}
iq_t = q_{xx} +2 [|q|^2 - \om^2] q +i \e [q_{xx} - \al q +\be ] \ ,
\label{spnls}
\end{equation}
where $q = q(t,x)$ is a complex-valued function of the two real 
variables $t$ and $x$, $t$ represents time, and $x$ represents
space. $q(t,x)$ is subject to periodic boundary condition of period 
$2 \pi$, and even constraint, i.e., 
\[
q(t,x + 2 \pi) = q(t,x)\ , \ \ q(t,-x) = q(t,x)\ .
\]
$\om \in (1/2, 1)$ is a positive constant, $\al >0$ and $\be >0$ 
are constants, and $\e > 0$ is the perturbation parameter. 
The following theorem was proved in \cite{Li01b}.
\begin{theorem}
There exists a $\e_0 > 0$, such that for any $\e \in (0, \e_0)$, there 
exists a codimension 1 surface in the external parameter space 
$(\alpha,\beta, \om) \in \RR^+\times  \RR^+\times 
\RR^+$ where $\om \in (\frac{1}{2}, 1)/S$, $S$ is a finite subset, and 
$\al \om < \be$. For any $(\alpha ,\beta, \omega)$ on the codimension 1
surface, the singularly perturbed nonlinear Schr\"odinger equation 
(\ref{spnls}) possesses a symmetric pair of Silnikov homoclinic orbits 
asymptotic 
to a saddle $Q_\epsilon$. The codimension 1 surface has the approximate 
representation given by $\al = 1/\k(\om)$, where $\k(\om)$ is plotted 
in Figure \ref{kappa}.
\label{shorbit}
\end{theorem}
The proof of the theorem is also completed 
through two measurements: the Melnikov measurement and the second measurement.
But more powerful machineries are needed \cite{Li01b}.

\subsection{Silnikov Homoclinic Orbits in Vector NLS Under Perturbations}

In recent years, novel results have been obtained on the solutions of 
the vector nonlinear Schr\"odinger equations \cite{AOT99} \cite{AOT00} 
\cite{YT01}. Abundant ordinary integrable results have been carried through 
\cite{WF00} \cite{FSW00}, including linear stability calculations 
\cite{FMMW00}. Specifically, the vector nonlinear 
Schr\"odinger equations can be written as
\begin{eqnarray*}
& & ip_t + p_{xx} + \frac{1}{2} (|p|^2 + \chi |q|^2) p = 0 , \\
& & iq_t + q_{xx} + \frac{1}{2} (\chi |p|^2 + |q|^2) q = 0 , 
\end{eqnarray*}
where $p$ and $q$ are complex valued functions of the two real variables 
$t$ and $x$, and $\chi$ is a positive constant. These equations describe
the evolution of two orthogonal pulse envelopes in birefringent optical 
fibers \cite{Men87} \cite{Men89}, with industrial applications in fiber 
communication systems \cite{HK95} and all-optical switching devices 
\cite{Isl92}. For linearly birefringent fibers \cite{Men87}, $\chi =2/3$.
For elliptically birefringent fibers, $\chi$ can take other positive 
values \cite{Men89}. When $\chi = 1$, these equations are first 
shown to be integrable by S. Manakov \cite{Man74}, and thus called Manakov 
equations. When $\chi$ is not 1 or 0, these equations are 
non-integrable. Propelled by the industrial applications, extensive 
mathematical studies on the vector nonlinear Schr\"odinger equations 
have been conducted. Like the scalar nonlinear Schr\"odinger equation, the 
vector nonlinear Schr\"odinger equations also possess figure eight 
structures in their phase space. Consider the singularly perturbed vector 
nonlinear Schr\"odinger equations,
\begin{eqnarray}
& & ip_t + p_{xx} + \frac{1}{2} [(|p|^2 + |q|^2)-\om^2] p = 
i \e [ p_{xx} -\al p - \be ]\ , \label{pnls1}\\
& & iq_t + q_{xx} + \frac{1}{2} [(|p|^2 + |q|^2)-\om^2] q = 
i \e [ q_{xx} -\al q - \be ]\  , \label{pnls2}
\end{eqnarray}
where $p(t,x)$ and $q(t,x)$ are subject to periodic boundary condition 
of period $2\pi$, and are even in $x$, i.e. 
\[
p(t,x + 2 \pi) = p(t,x)\ , \ \ p(t,-x) = p(t,x)\ , 
\]
\[
q(t,x + 2 \pi) = q(t,x)\ , \ \ q(t,-x) = q(t,x)\ , 
\]
$\om \in (1,2)$, $\al > 0$ and $\be$ are real constants, and 
$\e > 0$ is the perturbation parameter. We have
\begin{theorem}[\cite{Li02d}]
There exists a $\e_0 > 0$, such that 
for any $\e \in (0, \e_0)$, there exists 
a codimension 1 surface in the space of $(\alpha,\beta, \om) \in 
\RR^+\times \RR^+\times \RR^+$ where 
$\om \in (1, 2)/S$, $S$ is a finite subset, and 
$\al \om < \sqrt{2} \be$. For any $(\alpha ,\beta, \omega)$ on the 
codimension 1 surface, the singularly perturbed vector nonlinear 
Schr\"odinger equations (\ref{pnls1})-(\ref{pnls2}) possesses a 
homoclinic orbit asymptotic to a saddle
$Q_\epsilon$. This orbit is also the homoclinic orbit 
for the singularly perturbed scalar nonlinear Schr\"odinger equation
studied in last section, and is the only one asymptotic 
to the saddle $Q_\epsilon$ for the singularly perturbed 
vector nonlinear Schr\"odinger equations (\ref{pnls1})-(\ref{pnls2}).
The codimension 1 surface has the 
approximate representation given by $\al = 1/\k(\om)$, where $\k(\om)$ 
is plotted in Figure \ref{kappa}.
\end{theorem} 

\subsection{Silnikov Homoclinic Orbits in Discrete 
NLS Under Perturbations \label{hordnls}}

Consider the following perturbed discrete cubic nonlinear Schr\"odinger 
equations \cite{LM97}, 
\begin{eqnarray}
i\dot{q_n}&=&{1 \over h^2}\bigg[q_{n+1}-2q_n+q_{n-1}\bigg]+|q_n|^2(q_{n+1}+
           q_{n-1})-2\om^2 q_n \nonumber \\
        & &+i\e \bigg[-\al q_n +{1 \over h^2}(q_{n+1}-2q_n+q_{n-1})
           + \be \bigg], \label{PDNLS}
\end{eqnarray}
\nid
where $i=\sqrt{-1}$, $q_n$'s are complex variables,
\[
q_{n+N}=q_n, \ \ (\mbox{periodic}\  \mbox{condition}); \quad \mbox{and}\ 
q_{-n}=q_n, \ \ (\mbox{even}\  \mbox{condition});
\]
$h={1\over N}$, and 
\begin{eqnarray*}
& & N\tan{\pi \over N}< \om <N\tan{2\pi \over N},\ \ \mbox{for}\ N>3,\\
& & 3\tan{\pi \over 3}< \om < \infty, \ \ \mbox{for}\ N=3. \\
& & \e\in[0,\e_1),\ \alpha\ (>0), \ \be\ (>0) \ 
\mbox{are}\ \mbox{constants.}
\end{eqnarray*}
This is a $2(M+1)$ dimensional system, where
\[
M=N/2,\ \ (N\ \mbox{even}); \quad \mbox{and}\ 
M=(N-1)/2, \ \ (N\ \mbox{odd}).
\]
This system is a finite-difference discretization of the perturbed NLS 
(\ref{spnls}). The following theorem was proved in \cite{LM97}.

Denote by $\Sg_N\ (N\geq 7)$ the external parameter space,
\begin{eqnarray*}
\Sg_N&=&\bigg\{ (\om,\al,\be)\ \bigg | \ \om \in (N\tan{\pi \over N}, 
               N\tan{2\pi \over N}),\\ 
     & &\al\in (0,\al_0), \be\in (0,\be_0); \\
     & &\mbox{where}\ \al_0\ \mbox{and}\ \be_0\ \mbox{are}\ \mbox{any}
        \ \mbox{fixed}\ \mbox{positive}\ \mbox{numbers}. \bigg\}
\end{eqnarray*}
\begin{theorem}
For any $N$ ($7\leq N<\infty$), there exists a positive number $\e_0$, 
such that for any $\e \in (0,\e_0)$, there exists a 
codimension $1$ surface $E_\e$ in $\Sg_N$; for any
external parameters ($\om,\al,\be$) on $E_\e$, there exists a 
homoclinic orbit asymptotic to a saddle $Q_\e$.
The codimension $1$ surface $E_\e$ has the approximate expression
$\al=1/\k$, where $\k=\k(\om;N)$
is shown in Fig.\ref{nkappa}.
\label{dhorbit}
\end{theorem}
In the cases ($3\leq N \leq 6$), $\k$ is always negative. For 
$N \geq 7$, $\k$ 
can be positive as shown in Fig.\ref{nkappa}. When $N$ is even and 
$\geq 7$, there is in fact a symmetric pair of homoclinic orbits asymptotic to 
a fixed point $Q_\e$ at the same values of the external parameters; 
since for even $N$, we have the symmetry:
If $q_n=f(n,t)$ solves (\ref{PDNLS}),
then $q_n=f(n+N/2,t)$ also solves (\ref{PDNLS}). When $N$ is odd
and $\geq 7$, the study can not guarantee that two homoclinic orbits 
exist at the same value of the external parameters.
\begin{figure}
\caption{The graph of $\k(\om;N)$.}
\label{nkappa}
\end{figure}

\subsection{Comments on DSII Under Perturbations}

Consider the perturbed Davey-Stewartson II equations \cite{Li02b},
\begin{eqnarray*}
iq_t &=& \Upsilon q+ \bigg [2(|q|^2-\omega^2)+ u_y \bigg ]q
+i\epsilon f \ , \\
& & \ \ \Delta u = -4\partial_y |q|^2 \ ,
\end{eqnarray*}
where $q$ is a complex-valued function of the three variables ($t,x,y$), 
$u$ is a real-valued function of the three variables ($t,x,y$), 
$\Upsilon =\partial_{xx}-\partial_{yy}$, $\Delta=\partial_{xx}
+\partial_{yy}$, $\omega >0$ is a constant. We also consider periodic 
boundary conditions. The perturbation $f$ can be for example a singular 
perturbation 
\[
f = \Dl q - \al q + \be \ ,
\]
where $\alpha >0$, $\beta >0$ are constants, or a regular perturbation 
by mollifying $\Dl$ into a bounded Fourier multiplier
\[
\hat{\Dl} q = -\sum_{k \in Z^2} \be_k |k|^2 \tilde{q}_k \cos k_1 x 
\cos k_2 y \ ,
\]
in the case of periods ($2\pi , 2\pi$),
\[
\be_k = 1, \ \ |k| \leq N , \quad \be_k = |k|^{-2},\ \ |k| > N,
\]
for some large $N$, $|k|^2 = k_1^2 + k_2^2$.

Under both regular and singular 
perturbations, the rigorous Melnikov measurement can be established 
\cite{Li02b}. It turns out that only local well-posedness is necessary 
for rigorously setting up the Melnikov measurement, thanks to the fact 
that the unperturbed homoclinic orbits are 
classical solutions. The obstacle toward proving the existence of homoclinic 
orbits comes from a technical difficulty in solving a linear system 
to get the normal form for proving the size estimate of the stable 
manifold of a saddle. For details, see \cite{Li02b}.

\subsection{Transversal Homoclinic Orbits in a Periodically Perturbed 
Sine-Gordon Equation \label{PPSGE}}

Transversal homoclinic orbits in continuous systems often appear 
in two types of systems: (1). periodic systems where the Poincar\'e 
period map has a transversal homoclinic orbit, (2). autonomous 
systems where the homoclinic orbit is asymptotic to a hyperbolic 
limit cycle. 

Consider the periodically perturbed sine-Gordon (SG) equation,
\begin{equation}
u_{tt}=c^2 u_{xx}+\sin u+\epsilon [-a u_t+u^3 \chi(\| u\|)\cos t],
\label{PSG}
\end{equation} 
where
\[
\chi(\| u\|)=\left\{ \begin{array}{ll} 1, & \| 
u\| \leq M,\\ 0, & \| u\| \geq 2M,\end{array}\right.
\]
for $M<\| u\| <2M$, $\chi (\| u\|)$ is a smooth bump function, 
under odd periodic boundary condition,
\[
u(x+2\pi ,t)=u(x,t),\quad 
u(x,t)=-u(x,t),
\]
$\frac{1}{4}<c^2<1$, $a>0$, $\epsilon$ is a small 
perturbation parameter.
\begin{theorem}[\cite{LMSW96}, \cite{SZ00}] There exists an 
interval $I\subset \RR^{+}$ such that for any $a\in I$, there exists a
transversal homoclinic orbit
$u=\xi (x,t)$ asymptotic to $0$ in $H^{1}$.
\end{theorem}

\subsection{Transversal Homoclinic Orbits in a Derivative NLS \label{ADNS}}

Consider the derivative nonlinear Schr\"odinger equation,
\begin{equation}
i q_t = q_{xx} + 2 |q|^2 q +i \e \bigg [ (\frac{9}{16}-|q|^2 )q +\mu 
|\hat{\pa}_x q|^2 \bar{q} \bigg ]\ , \label{derNLS}
\end{equation}
where $q$ is a complex-valued function of two real variables $t$ and $x$,
$\e > 0$ is the perturbation parameter, $\mu$ is a real constant, and
$\hat{\pa}_x $ is a bounded Fourier multiplier,
\[
\hat{\pa}_x q = -\sum_{k=1}^K k \tq_k \sin kx\ , \quad 
\mbox{for} \ q = \sum_{k=0}^\infty \tq_k \cos kx\ ,
\]
and some fixed $K$. Periodic boundary condition 
and even constraint are imposed,
\[
q(t,x+2\pi ) = q(t,x)\ , \ \ q(t,-x)=q(t,x) \ . 
\] 
\begin{theorem}[\cite{Li02a}]
There exists a $\e_0 > 0$, such that 
for any $\e \in (0, \e_0)$, and $|\mu | > 5.8$,
there exist two transversal homoclinic orbits asymptotic to 
the limit cycle $q_c = \frac{3}{4} \exp \{ -i [ \frac{9}{8} t + \ga ]\}$.
\label{thmdns}
\end{theorem}

\section{Existence of Chaos}

The importance of homoclinic orbits with respect to chaotic dynamics was 
first realized by Poincar{\'e} \cite{Poi99}. In 1961, Smale constructed 
the well-known horseshoe in the neighborhood of a transversal homoclinic 
orbit \cite{Sma61} \cite{Sma65} \cite{Sma67}. In particular, 
Smale's theorem implies Birkhoff's theorem on the existence of a 
sequence of structurely stable periodic orbits in the neighborhood 
of a transversal homoclinic orbit \cite{Bir12}. In 1984 and 1988 \cite{Pal84}
\cite{Pal88}, Palmer gave a beautiful proof of Smale's theorem using a 
shadowing lemma. Later, this proof was generalized to infinite dimensions 
by Steinlein and Walther \cite{SW89} \cite{SW90} and Henry \cite{Hen94}. 
In 1967, Silnikov proved Smale's theorem for autonomous systems in finite 
dimensions using a fixed point argument \cite{Sil67b}. In 1996, Palmer proved 
Smale's theorem for autonomous systems in finite dimensions using shadowing 
lemma \cite{CKP95} \cite{Pal96}. In 2002, Li proved Smale's theorem for 
autonomous systems in infinite dimensions using shadowing lemma \cite{Li02a}.

For nontransversal homoclinic orbits, the most well-known type which leads 
to the existence of Smale horseshoes is the so-called Silnikov homoclinic 
orbit \cite{Sil65} \cite{Sil67a} \cite{Sil70} \cite{Den89} \cite{Den93}. 
Existence of Silnikov homoclinic orbits and new constructions of Smale 
horseshoes for concrete nonlinear wave systems have been established 
in finite dimensions \cite{LM97} \cite{LW97} and in infinite dimensions 
\cite{LMSW96} \cite{Li99a} \cite{Li01b}. 

\subsection{Shift Automorphism}

Let $\W$ be a set which consists of elements of the doubly infinite 
sequence form:
\[
a =(\cdot \cdot \cdot  a_{-2} a_{-1} a_0, a_1 a_2 \cdot \cdot \cdot ),
\]
where $a_k \in \{ l_1, \cdots, l_m\}$, m labels, and 
$k\in \ZZ$. We introduce a topology in $\W$
by taking as neighborhood basis of 
\[
a^* =( \cdot \cdot \cdot a^*_{-2} a^*_{-1} a^*_0, a^*_1 a^*_2 
\cdot \cdot \cdot ),
\]
the set 
\[
W_j =  \bigg \{ a\in \W \ \bigg | \ a_k=a^*_k\ (|k|<j) \bigg \}
\]
\nid
for $j=1,2,\cdot \cdot \cdot $. This makes $\W$ a topological space.
The shift automorphism $\chi$ is defined on $\W$ by
\begin{eqnarray*}
\chi &:& \W \mapsto \W, \\
  & & \forall a \in \W,\ \chi(a) = b,\ \mbox{where}\ b_k=a_{k+1}.
\end{eqnarray*}
The shift automorphism $\chi$ exhibits {\em{sensitive dependence on 
initial conditions}}, which is a hallmark of {\em{chaos}}.

\subsection{NLS Under Regular Perturbations
\label{shrnls}}

We continue from section \ref{horrnls} and consider the regularly perturbed 
nonlinear Schr\"odinger (NLS) equation (\ref{rpnls}). Starting from the 
Homoclinic Orbit Theorem \ref{rhorbit},
We have the following theorem on the existence of chaos. 
\begin{theorem}[Chaos Theorem, \cite{Li99a}]
Under certain generic assumptions for the 
perturbed nonlinear Schr\"odinger system (\ref{rpnls}), there exists a 
compact Cantor subset 
$\La$ of $H^1$ (the Sobolev space), $\La$ consists of points, and is 
invariant under a Poincar\'e map $P$.
$P$ restricted to $\La$, is topologically conjugate to the shift 
automorphism $\chi$ on four symbols $1, 2, -1, -2$. That is, there exists
a homeomorphism
\[
\phi \ : \ \W \mapsto  \La,
\]
\nid
such that the following diagram commutes:
\begin{equation} 
\begin{array}{ccc}
\W &\maprightu{\phi} & \Lambda\\
\mapdownl{\chi} & & \mapdownr{P}\\
\W & \maprightd{\phi} & \Lambda
\end{array} 
\nonumber
\end{equation}
\label{horserhthm}
\end{theorem}
Although the symmetric pair of Silnikov homoclinic orbits is not structurally 
stable, the Smale horseshoes are structurally stable. Thus, the Cantor 
sets and the conjugacy to the shift automorphism are structurally stable.

\subsection{NLS Under Singular Perturbations
\label{shsnls}}

We continue from section \ref{horsnls} and consider the singularly perturbed 
nonlinear Schr\"odinger (NLS) equation (\ref{spnls}). Starting from the 
Homoclinic Orbit Theorem \ref{shorbit}, we have
\begin{theorem}[Chaos Theorem, \cite{Li02c}]
Under certain generic assumptions for the 
perturbed nonlinear Schr\"odinger system (\ref{spnls}), there exists a 
compact Cantor subset 
$\La$ of $H^1$, $\La$ consists of points, and is invariant under a 
Poincar\'e map $P$.
$P$ restricted to $\La$, is topologically conjugate to the shift 
automorphism $\chi$ on four symbols $1, 2, -1, -2$. That is, there exists
a homeomorphism
\[
\phi \ : \ \W \mapsto  \La,
\]
\nid
such that the following diagram commutes:
\begin{equation} 
\begin{array}{ccc}
\W &\maprightu{\phi} & \Lambda\\
\mapdownl{\chi} & & \mapdownr{P}\\
\W & \maprightd{\phi} & \Lambda
\end{array} 
\nonumber
\end{equation}
\label{horseshthm}
\end{theorem}

\subsection{Discrete NLS Under Perturbations
\label{shdnls}}

We continue from section \ref{hordnls}. Starting from the 
Homoclinic Orbit Theorem \ref{dhorbit}, we have the theorem 
on the existence of chaos 
for the perturbed discrete nonlinear Schr\"odinger equation 
(\ref{PDNLS}). 
\begin{theorem}[Chaos Theorem, \cite{LW97}]
Under certain generic assumptions for the 
perturbed discrete nonlinear Schr\"odinger system (\ref{PDNLS}),
there exists a compact Cantor subset 
$\La$ of $\RR^{2(M+1)}$, $\La$ consists of points, and is invariant under 
a Poincar\'e map $P$.
$P$ restricted to $\La$, is topologically conjugate to the shift 
automorphism $\chi$ on two symbols $0, 1$. That is, there exists
a homeomorphism
\[
\phi \ : \ \W \mapsto  \La,
\]
\nid
such that the following diagram commutes:
\begin{equation} 
\begin{array}{ccc}
\W &\maprightu{\phi} & \Lambda\\
\mapdownl{\chi} & & \mapdownr{P}\\
\W & \maprightd{\phi} & \Lambda
\end{array} 
\nonumber
\end{equation}
\label{horsedhthm}
\end{theorem}

\subsection{Numerical Simulation of Chaos}

The finite-difference discretization of both the regularly and the singularly 
perturbed nonlinear Schr\"odinger equations (\ref{rpnls}) and (\ref{spnls})
leads to the same discrete perturbed nonlinear Schr\"odinger equation 
(\ref{PDNLS}).
 
In the chaotic regime, typical numerical output is shown in 
Figure \ref{numo}. Notice that there are two typical profiles at a fixed 
time: (1). a breather type profile with its hump located at the center of 
the spatial period interval, (2). a breather type profile with its hump 
located at the boundary (wing) of the spatial period interval. These two 
types of profiles are half spatial period translate of each other.
If we label the profiles, with their humps at the center of the spatial 
period interval, by ``C''; those profiles, with their humps at the wing
of the spatial period interval, by ``W'', then 
\begin{equation}
``W\mbox{''}=\sg \circ ``C\mbox{''}, \label{CWJ}
\end{equation}
where $\sg$ is the symmetry group element representing half spatial period 
translate. The time series of the output in Figure \ref{numo} is a 
chaotic jumping between ``C'' and ``W'', which we call ``chaotic center-wing 
jumping''. We interpret the chaotic center-wing jumping as the 
numerical realization of the shift automorphism $\chi$ on 
symbols. We can make this more precise in terms of the 
phase space geometry. The figure-eight structure of the integrable NLS 
projected 
onto the plane of the Fourier component $\cos x$ is illustrated in 
Figure \ref{pfi}, and labeled by $L_C$ and $L_W$. From 
the symmetry, we know that
\[
L_W = \sg \circ L_C.
\]
$L_C$ has the spatial-temporal profile realization
as in Figure \ref{figlc} with the hump located at the center. $L_W$ 
corresponds to the half spatial period translate 
of the spatial-temporal profile realization as in Figure \ref{figlc},
with the hump located at the boundary (wing).
An orbit inside $L_C$, $L_{Cin}$ has a spatial-temporal profile realization as
in Figure \ref{figlcin}. The half period translate of 
$L_{Cin}$, $L_{Win}$ is 
inside $L_W$. An orbit outside $L_C$ and $L_W$, $L_{out}$ has the
spatial-temporal profile realization as in Figure \ref{figlout}.
$S_l$ and $S_{l,\sg}$ are two phase blocks.
From these figures, we can see clearly that 
the chaotic center-wing jumping 
(Figure \ref{numo}) is the realization of the shift automorphism on 
symbols in $S_l \cup S_{l,\sg}$.
\begin{figure}
\caption{A chaotic solution in the discrete perturbed NLS 
system (\ref{PDNLS}).}
\label{numo}
\end{figure}
\begin{figure}
\caption{Figure-eight structure and its correspondence with the chaotic
center-wing jumping.}
\label{pfi}
\end{figure}
\begin{figure}
\caption{Spatial-temporal profile realization of $L_C$ in 
Figure \ref{pfi} (coordinates are the same as in Figure \ref{numo}).}
\label{figlc}
\end{figure}
\begin{figure}
\caption{Spatial-temporal profile realization of $L_{Cin}$ in 
Figure \ref{pfi} (coordinates are the same as in Figure \ref{numo}).}
\label{figlcin}
\end{figure}
\begin{figure}
\caption{Spatial-temporal profile realization of $L_{out}$ in 
Figure \ref{pfi} (coordinates are the same as in Figure \ref{numo}).}
\label{figlout}
\end{figure}

\subsection{Shadowing Lemma and Chaos in Finite-D Periodic 
Systems \label{Palmer}}

Since its invention \cite{Ano67}, shadowing lemma has been a useful tool
for solving many dynamical system problems. Here we only focus upon its 
use in proving the existence of chaos in a neighborhood of a transversal 
homoclinic orbit. According to the type of the system, the level of 
difficulty in proving the existence of chaos with a shadowing lemma
is different. 

A finite-dimensional periodic system can be written in the 
general form
\[
\dot{x} = F(x,t)\ ,
\]
where $x \in \RR^n$, and $F(x,t)$ is periodic in $t$. Let $f$ be the 
Poincar\'e period map.
\begin{definition}
A doubly infinite sequence $\{ y_j \}$ in $\RR^n$ is a $\dl$ pseudo-orbit 
of a $C^1$ diffeomorphism $f : \RR^n \mapsto \RR^n$ if for all integers $j$
\[
|y_{j+1}-f(y_j)| \leq \dl \ .
\]
An orbit $\{ f^j(x) \}$ is said to $\e$-shadow the $\dl$ pseudo-orbit 
$\{ y_j \}$ if for all integers $j$
\[
|f^j(x)-y_j| \leq \e \ .
\]
\end{definition}
\begin{definition}
A compact invariant set $S$ is hyperbolic if there are positive 
constants $K$, $\al$ and a projection matrix valued function $P(x)$, 
$x \in S$, of constant rank such that for all $x$ in $S$
\[
P(f(x))Df(x) = Df(x)P(x)\ ,
\]
\[
|Df^j(x)P(x)| \leq Ke^{-\al j}\ ,\quad (j \geq 0 )\ , 
\]
\[
|Df^j(x)(I-P(x))| \leq Ke^{\al j}\ ,\quad (j \leq 0 )\ . 
\] 
\end{definition}
\begin{theorem}[Shadowing Lemma \cite{Pal88}]
Let $S$ be a compact hyperbolic set for the $C^1$ diffeomorphism 
$f : \RR^n \mapsto \RR^n$. Then given $\e >0$ sufficiently small there 
exists $\dl >0$ such that every $\dl$ pseudo-orbit in $S$ has a 
unique $\e$-shadowing orbit.
\label{shal}
\end{theorem}
The proof of this theorem by Palmer \cite{Pal88} is overall a 
fixed point argument with the help of Green functions for linear maps.

Let $y_0$ be a transversal homoclinic point asymptotic to a saddle $x_0$ 
of a $C^1$ diffeomorphism $f : \RR^n \mapsto \RR^n$. Then the set
\[
S= \{ x_0 \} \cup \{ f^j(y_0): j \in Z\}
\]
is hyperbolic. Denote by $A_0$ and $A_1$ the two orbit segments of 
length $2m+1$
\[
A_0= \{ x_0,x_0, \cdots, x_0 \}\ , \quad 
A_1= \{ f^{-m}(y_0), f^{-m+1}(y_0), \cdots, f^{m-1}(y_0), f^{m}(y_0) \} \ .
\]
Let 
\[
a=(\cdots, a_{-1},a_0, a_1, \cdots ) \ , 
\]
where $a_j \in \{ 0,1\}$, be any doubly infinite binary sequence. Let 
$A$ be the doubly infinite sequence of points in $S$, associated with 
$a$
\[
A=\{ \cdots, A_{a_{-1}},A_{a_0}, A_{a_1}, \cdots \} \ .
\]
When $m$ is sufficiently large, $A$ is a $\dl$ pseudo-orbit in $S$.
By the shadowing lemma (Theorem \ref{shal}), there is a unique 
$\e$-shadowing orbit that shadows $A$. In this manner, Palmer \cite{Pal88}
gave a beautiful proof of Smale's horseshoe theorem.
\begin{definition}
Denote by $\Sg$ the set of doubly infinite binary sequences
\[
a=(\cdots, a_{-1},a_0, a_1, \cdots ) \ , 
\]
where $a_j \in \{ 0,1\}$. We give the set $\{ 0,1\}$ the discrete topology 
and $\Sg$ the product topology. The Bernoulli shift $\chi$ is defined by
\[
[\chi(a)]_j = a_{j+1}\ .
\]
\end{definition}
\begin{theorem}
Let $y_0$ be a transversal homoclinic point asymptotic to a saddle $x_0$ 
of a $C^1$ diffeomorphism $f : \RR^n \mapsto \RR^n$. Then there is a 
homeomorphism $\phi$ of $\Sg$ onto a compact subset of $\RR^n$ which is 
invariant under $f$ and such that when $m$ is sufficiently large
\[
f^{2m+1} \circ \phi = \phi \circ \chi\ ,
\]
that is, the action of $f^{2m+1}$ on $\phi(\Sg)$ is topologically conjugate 
to the action of $\chi$ on $\Sg$.
\end{theorem}
Here one can define $\phi(a)$ to be the point on the shadowing orbit that 
shadows the midpoint of the orbit segment $A_{a_0}$, which is either $x_0$ 
or $y_0$. The topological conjugacy can be easily verified. For details, 
see \cite{Pal88}.
Other references can be found in \cite{Pal84} \cite{Zen95}.

\subsection{Shadowing Lemma and Chaos in Infinite-D Periodic Systems}

An infinite-dimensional periodic system defined in a Banach space $X$ 
can be written in the general form
\[
\dot{x} = F(x,t)\ ,
\]
where $x \in X$, and $F(x,t)$ is periodic in $t$. Let $f$ be the 
Poincar\'e period map. When $f$ is a $C^1$ map which needs not to be 
invertible, shadowing lemma and symbolic dynamics around a transversal 
homoclinic orbit can both be established \cite{SW89} \cite{SW90} \cite{Hen94}.
Other references can be found in \cite{HL86} \cite{CLP89} \cite{Zen97} 
\cite{Bla86}. There exists also a work on horseshoe construction without
shadowing lemma for 
sinusoidally forced vibrations of buckled beam \cite{HM81}.

\subsection{Periodically Perturbed Sine-Gordon (SG) Equation}

We continue from section \ref{PPSGE}, and use the notations in  
section \ref{Palmer}. For the periodically perturbed sine-Gordon 
equation (\ref{PSG}), the Poincar\'e period map is a $C^1$ diffeomorphism 
in $H^1$. As a corollary of the result in last section, we have the theorem 
on the existence of chaos.
\begin{theorem} There is an integer $m$ and a homeomorphism $\phi$ 
of $\Sigma$ onto a compact Cantor subset $\Lambda $ of
$H^{1}$. $\Lambda$ is invariant under the Poincar\'e period-$2\pi$ map $P$ of 
the periodically perturbed sine-Gordon equation (\ref{PSG}). The 
action of $P^{2m+1}$ on $\Lambda$ is topologically conjugate to
the action of $\chi$ on
$\Sigma: P^{2m+1} \circ \phi =\phi \circ
\chi$. That is, the following diagram commutes:
\[
\begin{array}{ccc}
\Sg &\maprightu{\phi} & \Lambda\\
\mapdownl{\chi} & & \mapdownr{P^{2m+1}}\\
\Sg & \maprightd{\phi} & \Lambda
\end{array} 
\]
\end{theorem}

\subsection{Shadowing Lemma and Chaos in Finite-D Autonomous Systems}

A finite-dimensional autonomous system can be written in the 
general form
\[
\dot{x} = F(x)\ ,
\]
where $x \in \RR^n$. In this case, a transversal homoclinic orbit 
can be an orbit asymptotic to a normally hyperbolic limit cycle.
That is, it is an orbit in the intersection of the stable and 
unstable manifolds of a normally hyperbolic limit cycle. Instead of 
the Poincar\'e period map as for periodic system, one may want 
to introduce the so-called Poincar\'e return map which is a map induced by 
the flow on a codimension 1 section which is transversal to the limit 
cycle. Unfortunately, such a map is not even well-defined in the 
neighborhood of the homoclinic orbit. This poses a challenging 
difficulty in extending the arguments as in the case of a Poincar\'e 
period map. In 1996, Palmer \cite{Pal96} completed a proof of a 
shadowing lemma and existence of chaos using Newton's method. It will 
be difficult to extend this method to infinite dimensions, since it 
used heavily differentiations in time. 
Other references can be found in \cite{CKP95} \cite{CKP97} \cite{Sil67b} 
\cite{FS77}.

\subsection{Shadowing Lemma and Chaos in Infinite-D Autonomous Systems}

An infinite-dimensional autonomous system defined in a Banach space $X$ 
can be written in the general form
\[
\dot{x} = F(x)\ ,
\]
where $x \in X$. In 2002, the author \cite{Li02a} completed a proof of a 
shadowing lemma and existence of chaos using Fenichel's persistence of 
normally hyperbolic invariant manifold idea. The setup is as follows,
\begin{itemize}
\item {\bf Assumption (A1):} There exist a hyperbolic limit cycle
$S$ and a transversal homoclinic orbit $\xi$ 
asymptotic to $S$. As curves, $S$ and $\xi$ are $C^{3}$.
\item {\bf Assumption (A2):} The Fenichel fiber theorem is valid at $S$. 
That is, there exist a family of unstable Fenichel fibers
$\{ {\mathcal F}^{u}(q): \  q \in S\}$ and a family of stable Fenichel 
fibers $\{ {\mathcal F}^{s}(q): \  q\in S\}$. For each fixed $q\in S$,
${\mathcal F}^{u}(q)$ and ${\mathcal F}^{s}(q)$ are $C^{3}$ submanifolds.
${\mathcal F}^{u}(q)$ and ${\mathcal F}^{s}(q)$ are $C^{2}$ in $q,
\forall q\in S$. The unions $\bigcup_{q\in S}{\mathcal F}^{u}(q)$ and
$\bigcup_{q\in S}{\mathcal F}^{s}(q)$ are the unstable and stable 
manifolds of $S$. Both families are invariant, i.e.
\[
F^{t}({\mathcal F}^{u}(q))\subset
{\mathcal F}^{u}(F^{t}(q)),  \forall\ t \leq 0, q\in S,
\]
\[
F^{t}({\mathcal F}^{s}(q))\subset {\mathcal F}^{s}(F^{t}(q)), 
\forall\ t \geq 0,  q \in S, 
\]
where $F^{t}$ is the evolution operator. There are positive constants
$\k$ and $\widehat{C}$ such that $\forall q\in S$, $\forall 
q^{-}\in {\mathcal F}^{u}(q)$ and $\forall q^{+}\in
{\mathcal F}^{s}(q)$,
\[
\| F^{t}(q^{-})-F^{t}(q)\| \leq
\widehat{C}e^{\k t}\| q^{-}-q\|,  \forall \ t \leq 0\ ,
\]
\[
\| F^{t}(q^{+})-F^{t}(q)\| \leq \widehat{C}e^{-\k t}\| q^{+}-q\|, 
\forall \ t \geq 0\ .
\]
\item {\bf Assumption (A3):} $F^{t}(q)$ is $C^{0}$ in $t$, for
$t\in (-\infty ,\infty)$, $q\in X$. For any fixed $t\in 
(-\infty ,\infty )$, $F^{t}(q)$ is a $C^{2}$ diffeomorphism on
$X$.
\end{itemize}
\begin{remark}
Notice that we do not assume that as functions of time, $S$ and $\xi$ 
are $C^3$ , and we only assume that as curves, $S$ and $\xi$ 
are $C^3$.
\end{remark}

Under the above setup, a shadowing lemma and existence of chaos can be
proved \cite{Li02a}. Another crucial element in the argument is the
establishment of a $\la$-lemma (also called inclination lemma) \cite{Li02a}. 

\subsection{A Derivative Nonlinear Schr\"odinger Equation}

We continue from section \ref{ADNS}, and consider the derivative nonlinear 
Schr\"odinger equation (\ref{derNLS}). The transversal homoclinic orbit 
given in Theorem \ref{thmdns} is a classical solution. Thus, Assumption (A1) 
is valid. Assumption (A2) follows from the standard arguments 
in \cite{LW97} \cite{LMSW96} \cite{Li01b}. Since the perturbation 
in (\ref{derNLS}) is bounded, Assumption (A3) follows from standard 
arguments. Thus there exists chaos in the derivative nonlinear 
Schr\"odinger equation (\ref{derNLS}) \cite{Li02a}.

\section{Lax Pairs of Euler Equations of Inviscid Fluids}

The governing equations for the incompressible viscous fluid flow are the 
Navier-Stokes equations. Turbulence occurs in the regime of high Reynolds 
number. By formally setting the Reynolds number equal to infinity, the 
Navier-Stokes equations reduce to the Euler equations of incompressible 
inviscid fluid flow. One may view the Navier-Stokes equations with large 
Reynolds number as a singular perturbation of the Euler equations.

Results of T. Kato show that 2D Navier-Stokes equations are globally 
well-posed in $C^0([0, \infty); H^s(R^2)), \ s>2$, and for any 
$0 < T < \infty$, the mild solutions of the 2D Navier-Stokes equations
approach those of the 2D Euler equations in $C^0([0, T]; H^s(R^2))$ 
\cite{Kat86}. 3D Navier-Stokes equations are locally well-posed in 
$C^0([0, \tau]; H^s(R^3)), \ s>5/2$, and the mild solutions of the 3D 
Navier-Stokes equations approach those of the 3D Euler equations in 
$C^0([0, \tau]; H^s(R^3))$, where $\tau$ depends on the norms of the initial 
data and the external force \cite{Kat72} \cite{Kat75}. Extensive studies on 
the inviscid limit have been carried by J. Wu et al. \cite{Wu96} 
\cite{CW96} \cite{Wu98} \cite{BW99}. There is no doubt that mathematical 
study on Navier-Stokes (Euler) equations is one of the most important 
mathematical problems. In fact, Clay Mathematics Institute has posted the 
global well-posedness of 3D Navier-Stokes equations as one of the one 
million dollars problems.

V. Arnold \cite{Arn66} realized that 2D Euler equations are a Hamiltonian 
system. Recently, the author found Lax pair structures for Euler 
equations \cite{Li01a} \cite{LY01} \cite{Li02e} \cite{Li02f} \cite{LS02}.
Understanding the structures of solutions to Euler equations is of 
fundamental interest. Of particular interest is the question on the 
global well-posedness of 3D Navier-Stokes and Euler equations. Our number 
one hope is that the Lax pair structures can be useful in investigating 
the global well-posedness. Our secondary hope is that the Darboux 
transformation \cite{LY01} associated with the Lax pair can generate 
explicit representation of homoclinic structures \cite{Li00a}. 

The philosophical significance of the existence of Lax pairs for Euler 
equations is even more important. If one defines integrability of an equation 
by the existence of a Lax pair, then both 2D and 3D Euler equations 
are integrable. More importantly, both 2D and 3D Navier-Stokes equations 
at high Reynolds numbers are singularly perturbed integrable systems. 
Such a point of view changes our old ideology on Euler and Navier-Stokes 
equations.

\subsection{A Lax Pair for 2D Euler Equation}

The 2D Euler equation can be written in the vorticity form,
\begin{equation}
\pa_t \Om + \{ \Psi, \Om \} = 0 \ ,
\label{euler}
\end{equation}
where the bracket $\{\ ,\ \}$ is defined as
\[
\{ f, g\} = (\pa_x f) (\pa_y g) - (\pa_y f) (\pa_x g) \ ,
\]
$\Om$ is the vorticity, and $\Psi$ is the stream function given by,
\[
u=- \pa_y \Psi \ ,\ \ \ v=\pa_x \Psi \ ,
\]
and the relation between vorticity $\Om$ and stream 
function $\Psi$ is,
\[
\Om =\pa_x v - \pa_y u =\Dl \Psi \ .
\]
\begin{theorem}[Li, \cite{Li01a}]
The Lax pair of the 2D Euler equation (\ref{euler}) is given as
\begin{equation}
\left \{ \begin{array}{l} 
L \varphi = \la \varphi \ ,
\\
\pa_t \varphi + A \varphi = 0 \ ,
\end{array} \right.
\label{laxpair}
\end{equation}
where
\[
L \varphi = \{ \Om, \varphi \}\ , \ \ \ A \varphi = \{ \Psi, \varphi \}\ ,
\]
and $\la$ is a complex constant, and $\varphi$ is a complex-valued function.
\label{2dlp}
\end{theorem}

\subsection{A Darboux Transformation for 2D Euler Equation}

Consider the Lax pair (\ref{laxpair}) at $\la =0$, i.e.
\begin{eqnarray}
& & \{ \Om, p \} = 0 \ , \label{d1} \\
& & \pa_t p + \{ \Psi, p \} = 0 \ , \label{d2} 
\end{eqnarray}
where we replaced the notation $\varphi$ by $p$.
\begin{theorem}[\cite{LY01}]
Let $f = f(t,x,y)$ be any fixed solution to the system 
(\ref{d1}, \ref{d2}), we define the Gauge transform $G_f$:
\begin{equation}
\tilde{p} = G_f p = \frac {1}{\Om_x}[p_x -(\pa_x \ln f)p]\ ,
\label{gauge}
\end{equation}
and the transforms of the potentials $\Om$ and $\Psi$:
\begin{equation}
\tilde{\Psi} = \Psi + F\ , \ \ \ \tilde{\Om} = \Om + \Dl F \ ,
\label{ptl}
\end{equation}
where $F$ is subject to the constraints
\begin{equation}
\{ \Om, \Dl F \} = 0 \ , \ \ \ \{ \Dl F, F \} = 0\ .
\label{constraint}
\end{equation}
Then $\tilde{p}$ solves the system (\ref{d1}, \ref{d2}) at 
$(\tilde{\Om}, \tilde{\Psi})$. Thus (\ref{gauge}) and 
(\ref{ptl}) form the Darboux transformation for the 2D 
Euler equation (\ref{euler}) and its Lax pair (\ref{d1}, \ref{d2}).
\label{dt}
\end{theorem}
\begin{remark}
For KdV equation and many other soliton equations, the 
Gauge transform is of the form \cite{MS91},
\[
\tilde{p} =  p_x -(\pa_x \ln f)p \ .
\]
In general, Gauge transform does not involve potentials.
For 2D Euler equation, a potential factor $\frac {1}{\Om_x}$
is needed. From (\ref{d1}), one has
\[
\frac{p_x}{\Om_x} = \frac{p_y}{\Om_y} \ .
\]
The Gauge transform (\ref{gauge}) can be rewritten as
\[
\tilde{p} = \frac{p_x}{\Om_x} - \frac{f_x}{\Om_x} \frac{p}{f}
=\frac{p_y}{\Om_y} - \frac{f_y}{\Om_y} \frac{p}{f}\ .
\]
The Lax pair (\ref{d1}, \ref{d2}) has a symmetry, i.e. it is 
invariant under the transform $(t,x,y) \ra (-t,y,x)$. The form 
of the Gauge transform (\ref{gauge}) resulted from the inclusion 
of the potential factor $\frac {1}{\Om_x}$, is consistent with 
this symmetry.
\end{remark}
Our hope is to use the Darboux transformation to generate homoclinic 
structures for 2D Euler equation \cite{Li00a}.

\subsection{A Lax Pair for Rossby Wave Equation}

The Rossby wave equation is
\[
\pa_t \Om + \{ \Psi , \Om \} + \be \pa_x \Psi = 0 \ ,
\]
where $\Om = \Om (t,x,y)$ is the vorticity, 
$\{ \Psi , \Om \} = \Psi_x \Om_y - \Psi_y \Om_x $, 
and $\Psi = \Dl^{-1} \Om$ 
is the stream function. Its Lax pair can be obtained
by formally conducting the transformation, $\Om = \tilde{\Om} +\be y$,
to the 2D Euler equation \cite{Li01a},
\[
\{ \Om , \varphi \} - \be \pa_x \varphi = \la \varphi \ ,
\quad \pa_t \varphi + \{ \Psi , \varphi \} = 0 \ ,
\]
where $\varphi$ is a complex-valued function, and $\la$ is 
a complex parameter. 

\subsection{Lax Pairs for 3D Euler Equation}

The 3D Euler equation can be written in vorticity form,
\begin{equation}
\pa_t \Om + (u \cdot \na) \Om - (\Om \cdot \na) u = 0 \ ,
\label{3deuler}
\end{equation}
where $u = (u_1, u_2, u_3)$ is the velocity, $\Om = (\Om_1, \Om_2, \Om_3)$
is the vorticity, $\na = (\pa_x, \pa_y, \pa_z)$, 
$\Om = \na \times u$, and $\na \cdot u = 0$. $u$ can be 
represented by $\Om$ for example through Biot-Savart law.
\begin{theorem}
The Lax pair of the 3D Euler equation (\ref{3deuler}) is given as
\begin{equation}
\left \{ \begin{array}{l} 
L \phi = \la \phi \ ,
\\
\pa_t \phi + A \phi = 0 \ ,
\end{array} \right.
\label{alaxpair}
\end{equation}
where
\[
L \phi = (\Om \cdot \na )\phi \ , 
\ \ \ A \varphi = (u \cdot \na )\phi \ , 
\]
$\la$ is a complex constant, and $\phi$ is a complex scalar-valued function.
\end{theorem}
\begin{theorem}[\cite{Chi00}]
Another Lax pair of the 3D Euler equation (\ref{3deuler}) is given as
\begin{equation}
\left \{ \begin{array}{l} 
L \varphi = \la \varphi \ ,
\\
\pa_t \varphi + A \varphi = 0 \ ,
\end{array} \right.
\label{3dlaxpair}
\end{equation}
where
\[
L \varphi = (\Om \cdot \na )\varphi - (\varphi \cdot \na )\Om \ , 
\ \ \ A \varphi = (u \cdot \na )\varphi - (\varphi \cdot \na ) u \ , 
\]
$\la$ is a complex constant, and $\varphi = (\varphi_1, \varphi_2, 
\varphi_3)$ is a complex 3-vector valued function.
\end{theorem}
Our hope is that the infinitely many conservation laws generated by $\la 
\in \CC$ can provide a priori estimates for the global well-posedness of 
3D Navier-Stokes equations, or better understanding on the global 
well-posedness \cite{LS02}. For more informations on the topics, 
see \cite{LY01} \cite{LS02}.

\section{Dynamical System Studies on 2D Euler Equation}

\subsection{Linearized 2D Euler Equation at a Fixed Point}

To begin an infinite dimensional dynamical system study, we investigate 
the linearized 2D Euler equation at a fixed point \cite{Li00b}
\cite{Li02e} \cite{Li02e}. We consider the 
2D Euler equation (\ref{euler}) under periodic boundary condition in 
both $x$ and $y$ directions with period $2\pi$. Expanding $\Om$ into 
Fourier series,
\[
\Om =\sum_{k\in \Z} \om_k \ e^{ik\cdot X}\ ,
\]
where $\om_{-k}=\overline{\om_k}\ $, $k=(k_1,k_2)^T$, 
and $X=(x,y)^T$. The 2D Euler equation
can be rewritten as 
\begin{equation}
\dot{\om}_k = \sum_{k=p+q} A(p,q) \ \om_p \om_q \ ,
\label{Keuler}
\end{equation}
where $A(p,q)$ is given by,
\begin{eqnarray}
A(p,q)&=& {1\over 2}[|q|^{-2}-|p|^{-2}](p_1 q_2 -p_2 q_1) \nonumber \\
\label{Af} \\      
      &=& {1\over 2}[|q|^{-2}-|p|^{-2}]\left | \begin{array}{lr} 
p_1 & q_1 \\ p_2 & q_2 \\ \end{array} \right | \ , \nonumber
\end{eqnarray}
where $|q|^2 =q_1^2 +q_2^2$ for $q=(q_1,q_2)^T$, similarly for $p$.

\begin{figure}[ht]
  \begin{center}
    \leavevmode
      \setlength{\unitlength}{2ex}
  \begin{picture}(36,27.8)(-18,-12)
    \thinlines
\multiput(-12,-11.5)(2,0){13}{\line(0,1){23}}
\multiput(-16,-10)(0,2){11}{\line(1,0){32}}
    \thicklines
\put(0,-14){\vector(0,1){28}}
\put(-18,0){\vector(1,0){36}}
\put(0,15){\makebox(0,0){$k_2$}}
\put(18.5,0){\makebox(0,0)[l]{$k_1$}}
\qbezier(-5.5,0)(-5.275,5.275)(0,5.5)
\qbezier(0,5.5)(5.275,5.275)(5.5,0)
\qbezier(5.5,0)(5.275,-5.275)(0,-5.5)
\qbezier(0,-5.5)(-5.275,-5.275)(-5.5,0)
    \thinlines
\put(4,4){\circle*{0.5}}
\put(0,0){\vector(1,1){3.7}}
\put(4.35,4.35){$p$}
\put(4,-4){\circle*{0.5}}
\put(8,0){\circle*{0.5}}
\put(-8,0){\circle*{0.5}}
\put(-8,-2){\circle*{0.5}}
\put(-12,-4){\circle*{0.5}}
\put(-12,-6){\circle*{0.5}}
\put(-4,2){\circle*{0.5}}
\put(-4,4){\circle*{0.5}}
\put(0,6){\circle*{0.5}}
\put(0,8){\circle*{0.5}}
\put(4,10){\circle*{0.5}}
\put(12,4){\circle*{0.5}}
\put(0,-8){\circle*{0.5}}
\put(-4,-12){\line(1,1){17.5}}
\put(-13.5,-7.5){\line(1,1){19.5}}
\put(-13.5,-5.5){\line(1,1){17.5}}
\put(-3.6,1.3){$\hat{k}$}
\put(-7,12.1){\makebox(0,0)[b]{$(-p_2, p_1)^T$}}
\put(-6.7,12){\vector(1,-3){2.55}}
\put(6.5,13.6){\makebox(0,0)[l]{$\Sg_{\hat{k}}$}}
\put(6.4,13.5){\vector(-2,-3){2.0}}
\put(7,-12.1){\makebox(0,0)[t]{$(p_2, -p_1)^T$}}
\put(6.7,-12.25){\vector(-1,3){2.62}}
\put(-4.4,-13.6){\makebox(0,0)[r]{$\bar{D}_{|p|}$}}
\put(-4.85,-12.55){\vector(1,3){2.45}}
\end{picture}
\end{center}
\caption{An illustration of the classes $\Sg_{\hk}$ and the disk 
$\bar{D}_{|p|}$.}
\label{class}
\end{figure}
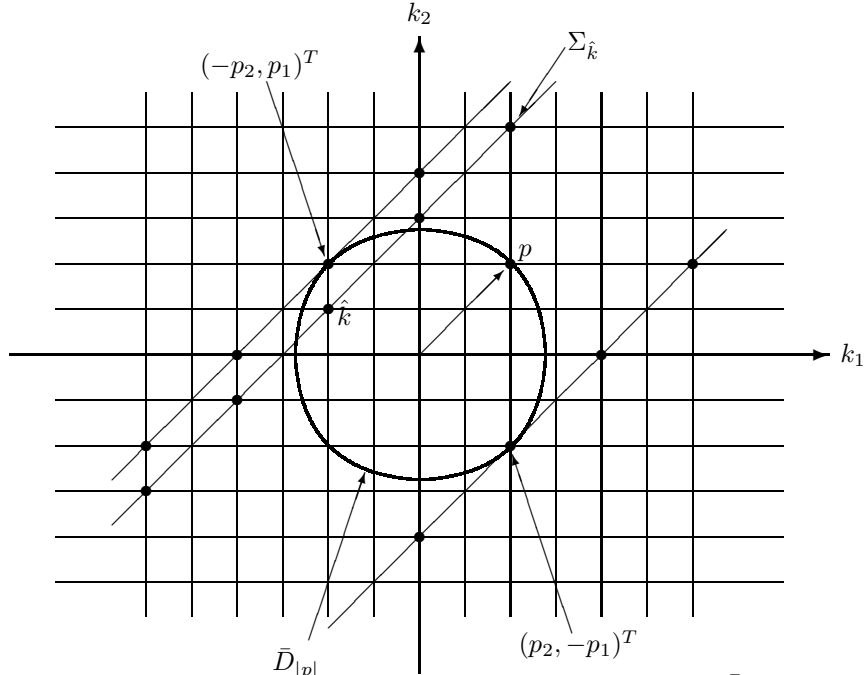

\begin{figure}[ht]
  \begin{center}
    \leavevmode
      \setlength{\unitlength}{2ex}
  \begin{picture}(36,27.8)(-18,-12)
    \thicklines
\put(0,-14){\vector(0,1){28}}
\put(-18,0){\vector(1,0){36}}
\put(0,15){\makebox(0,0){$\Im \{ \la \}$}}
\put(18.5,0){\makebox(0,0)[l]{$\Re \{ \la \}$}}
\put(0.1,-7){\line(0,1){14}}
\put(.2,-.2){\makebox(0,0)[tl]{$0$}}
\put(-0.2,-7){\line(1,0){0.4}}
\put(-0.2,7){\line(1,0){0.4}}
\put(2.0,-6.4){\makebox(0,0)[t]{$-i2|b|$}}
\put(2.0,7.6){\makebox(0,0)[t]{$i2|b|$}}
\put(2.4,3.5){\circle*{0.5}}
\put(-2.4,3.5){\circle*{0.5}}
\put(2.4,-3.5){\circle*{0.5}}
\put(-2.4,-3.5){\circle*{0.5}}
\put(5,4){\circle*{0.5}}
\put(-5,4){\circle*{0.5}}
\put(5,-4){\circle*{0.5}}
\put(-5,-4){\circle*{0.5}}
\put(8,6){\circle*{0.5}}
\put(-8,6){\circle*{0.5}}
\put(8,-6){\circle*{0.5}}
\put(-8,-6){\circle*{0.5}}
\end{picture}
  \end{center}
\caption{The spectrum of $\LL_{\hk}$ in case (2).}
\label{spla2}
\end{figure}
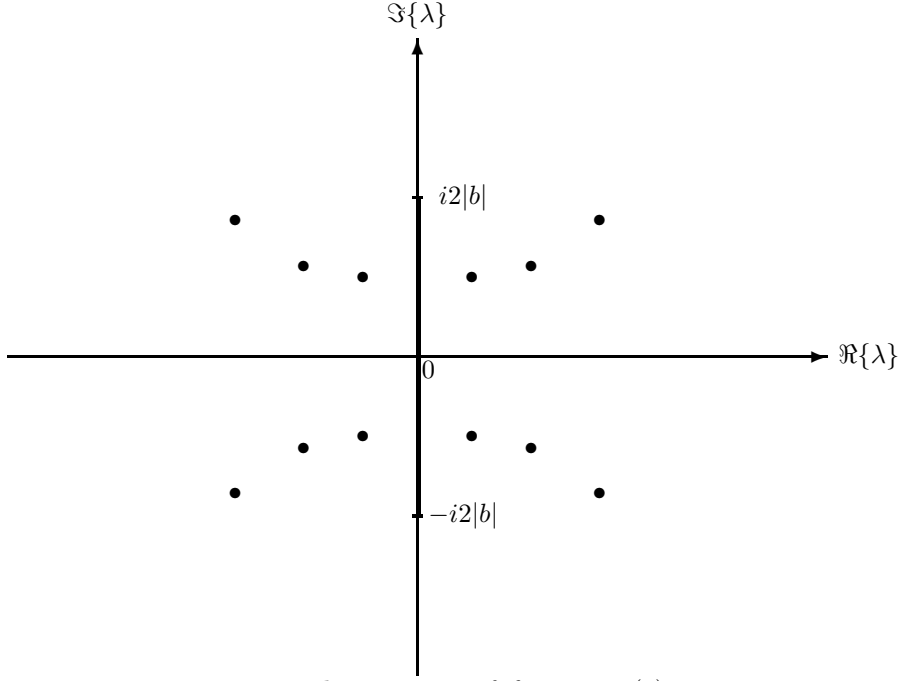

Denote $\{ \om_k \}_{k\in \Z}$ by $\om$. We consider the simple fixed point 
$\om^*$ \cite{Li00b}:
\begin{equation}
\om^*_p = \Ga,\ \ \ \om^*_k = 0 ,\ \mbox{if} \ k \neq p \ \mbox{or}\ -p,
\label{fixpt}
\end{equation}
of the 2D Euler equation (\ref{Keuler}), where 
$\Ga$ is an arbitrary complex constant. 
The {\em{linearized two-dimensional Euler equation}} at $\om^*$ is given by,
\begin{equation}
\dot{\om}_k = A(p,k-p)\ \Ga \ \om_{k-p} + A(-p,k+p)\ \bar{\Ga}\ \om_{k+p}\ .
\label{LE}
\end{equation}
\begin{definition}[Classes]
For any $\hk \in \Z$, we define the class $\Sg_{\hk}$ to be the subset of 
$\Z$:
\[
\Sg_{\hk} = \bigg \{ \hk + n p \in \Z \ \bigg | \ n \in \ZZ, \ 
\ p \ \mbox{is specified in (\ref{fixpt})} \bigg \}.
\]
\label{classify}
\end{definition}
\nid
See Fig.\ref{class} for an illustration of the classes. 
According to the classification 
defined in Definition \ref{classify}, the linearized two-dimensional Euler 
equation (\ref{LE}) decouples into infinitely many {\em{invariant subsystems}}:
\begin{eqnarray}
\dot{\omega}_{\hat{k} + np} &=& A(p, \hat{k} + (n-1) p) 
     \ \Gamma \ \omega_{\hat{k} + (n-1) p} \nonumber \\
& & + \ A(-p, \hat{k} + (n+1)p)\ 
     \bar{\Gamma} \ \omega_{\hat{k} +(n+1)p}\ . \label{CLE}
\end{eqnarray}
Let $\LL_{\hk}$ be the linear operator defined by the right hand side of 
(\ref{CLE}), and $H^s$ be the Sobolev space where $s \geq 0$ is an integer and 
$H^0=\ell_2$. 
\begin{theorem}[\cite{Li00b} \cite{Li02m}]
The eigenvalues of the linear operator $\LL_{\hk}$ in $H^s$ are of 
four types: real pairs ($c, -c$), purely imaginary pairs ($id, -id$), 
quadruples ($\pm c \pm id$), and zero eigenvalues.
\end{theorem}
\begin{definition}[The Disk]
The open disk of radius $|p|$ in $\Z$, denoted by $D_{|p|}$, is defined as
\[
  D_{|p|} = \bigg \{ k \in \Z \ \bigg| 
      \ |k| < |p| \bigg \} \, ,
\]
and the closure of $D_{|p|}$, denoted by
$\bar{D}_{|p|}$, is defined as
\[ 
 \bar{D}_{|p|} = \bigg \{ k \in \Z \ \bigg| 
     \ |k| \leq |p| \bigg \} \, .
\]
\end{definition}
\begin{theorem}[The Spectral Theorem, \cite{Li00b} \cite{Li02m}] We have 
the following claims on the spectrum of the linear operator $\LL_{\hk}$:
\begin{enumerate}
\item If $\Sg_{\hat{k}} \cap \bar{D}_{|p|} = \emptyset$, then the entire
$H^s$ spectrum of the linear operator $\LL_{\hk}$ 
is its continuous spectrum. See Figure \ref{splb}, where
$b= - \frac{1}{2}|\Gamma | |p|^{-2} 
\left|
  \begin{array}{cc}
p_1 & \hat{k}_1 \\
p_2 & \hat{k}_2
  \end{array}
\right| \ .$
That is, both the residual and the point spectra of $\LL_{\hk}$ are empty.
\item If $\Sg_{\hat{k}} \cap \bar{D}_{|p|} \neq \emptyset$, then the entire
essential $H^s$ spectrum of the linear operator $\LL_{\hk}$ is its 
continuous spectrum. 
That is, the residual 
spectrum of $\LL_{\hk}$ is empty. The point 
spectrum of $\LL_{\hk}$ is symmetric with respect to both real and 
imaginary axes. 
See Figure \ref{spla2}.
\end{enumerate}
\label{SST}
\end{theorem}
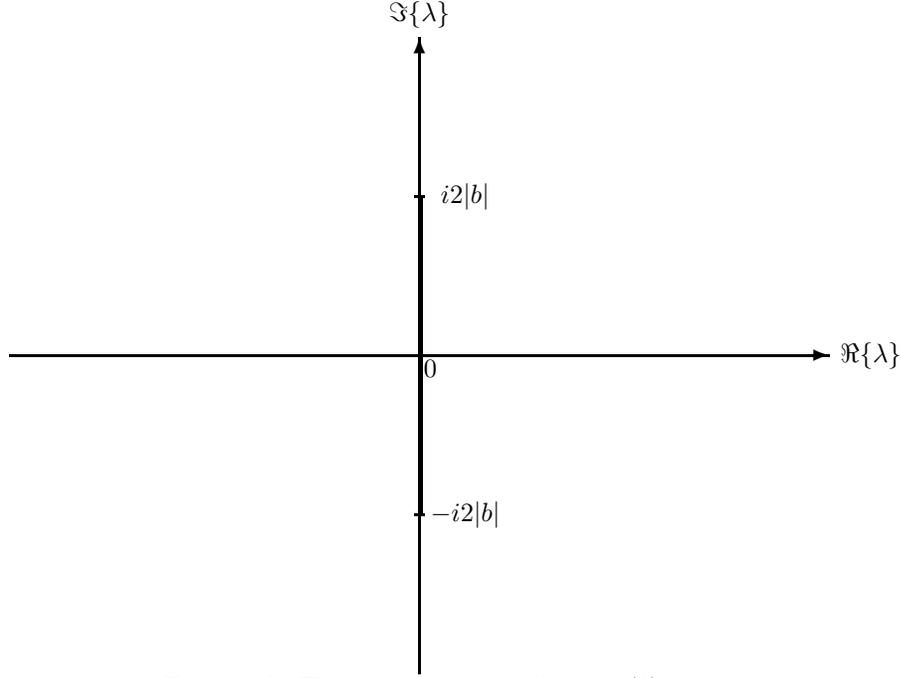
\begin{figure}[ht]
  \begin{center}
    \leavevmode
      \setlength{\unitlength}{2ex}
  \begin{picture}(36,27.8)(-18,-12)
    \thicklines
\put(0,-14){\vector(0,1){28}}
\put(-18,0){\vector(1,0){36}}
\put(0,15){\makebox(0,0){$\Im \{ \la \}$}}
\put(18.5,0){\makebox(0,0)[l]{$\Re \{ \la \}$}}
\put(0.1,-7){\line(0,1){14}}
\put(.2,-.2){\makebox(0,0)[tl]{$0$}}
\put(-0.2,-7){\line(1,0){0.4}}
\put(-0.2,7){\line(1,0){0.4}}
\put(2.0,-6.4){\makebox(0,0)[t]{$-i2|b|$}}
\put(2.0,7.6){\makebox(0,0)[t]{$i2|b|$}}
\end{picture}
  \end{center}
\caption{The spectrum of $\LL_{\hk}$ in case (1).}
\label{splb}
\end{figure}
\nid
Denote by $L$ the right hand side of (\ref{LE}), i.e. the whole 
linearized 2D Euler operator, the spectral mapping theorem 
holds.
\begin{theorem}[\cite{LLM01}]
$$\sigma(e^{tL})=e^{t\sigma(L)}, t\neq 0.$$
\end{theorem}
\nid
Let $\z$ denote the number of points $q \in
\Z$ that belong to the open disk of radius $|p|$, $D_{|p|}$, 
and such that $q$ is not parallel to $p$. 
\begin{theorem}[\cite{LLM01}]
The number of nonimaginary eigenvalues of $L$ (counting the multiplicities) 
does not exceed $2\z$.
\end{theorem}
\nid
Another interesting discussion upon the discrete spectrum can be found in 
\cite{Fad71}.

Since the introduction of continued fractions for calculating 
the eigenvalues of steady fluid flow, by Meshalkin and Sinai \cite{MS61},
this topics had been extensively explored \cite{Yud65} \cite{Liu92a} 
\cite{Liu92b} \cite{Liu93} \cite{Liu94a} \cite{Liu95} \cite{BFY99} 
\cite{Li00b}. Rigorous justification on the continued fraction calculation 
was given in \cite{Li00b} \cite{Liu95}.
As an example, we take $p=(1,1)^T$. When $\Ga \neq 0$, the fixed point 
has $4$ eigenvalues which form a 
quadruple. These four eigenvalues appear in the invariant 
linear subsystem labeled by $\hk = (-3,-2)^T$. One of them is \cite{Li00b}:
\begin{equation}
\tla=2 \lambda / | \Gamma | = 0.24822302478255 \ + \ i \ 0.35172076526520\ .
\label{evun}
\end{equation}
See Figure \ref{figev} for an illustration. The essential spectrum 
(= continuous spectrum) of $\LL_{\hk}$ with $\hk = (-3,-2)^T$ is the segment 
on the imaginary axis shown in Figure \ref{figev}, where 
$b = -\frac{1}{4} \Ga$. The essential spectrum 
(= continuous spectrum) of the linear 2D Euler operator at this fixed point 
is the entire imaginary axis. 
\begin{figure}[ht]
  \begin{center}
    \leavevmode
      \setlength{\unitlength}{2ex}
  \begin{picture}(36,27.8)(-18,-12)
    \thicklines
\put(0,-14){\vector(0,1){28}}
\put(-18,0){\vector(1,0){36}}
\put(0,15){\makebox(0,0){$\Im \{ \la \}$}}
\put(18.5,0){\makebox(0,0)[l]{$\Re \{ \la \}$}}
\put(2.4,3.5){\circle*{0.5}}
\put(-2.4,3.5){\circle*{0.5}}
\put(2.4,-3.5){\circle*{0.5}}
\put(-2.4,-3.5){\circle*{0.5}}  
\put(0.1,-10){\line(0,1){20}}
\put(.2,-.2){\makebox(0,0)[tl]{$0$}}
\put(-0.2,-10){\line(1,0){0.4}}
\put(-0.2,10){\line(1,0){0.4}}
\put(2.0,-9.4){\makebox(0,0)[t]{$-i2|b|$}}
\put(2.0,10.6){\makebox(0,0)[t]{$i2|b|$}}
\end{picture}
  \end{center}
\caption{The spectrum of $\LL_{\hk}$ with $\hk = (-3,-2)^T$, when $p=(1,1)^T$.}
\label{figev}
\end{figure}
A rather well-known open problem is proving the existence of unstable, stable,
and center manifolds. The main difficulty comes from the fact that the 
nonlinear term is non-Lipschitzian.

\subsection{Models}

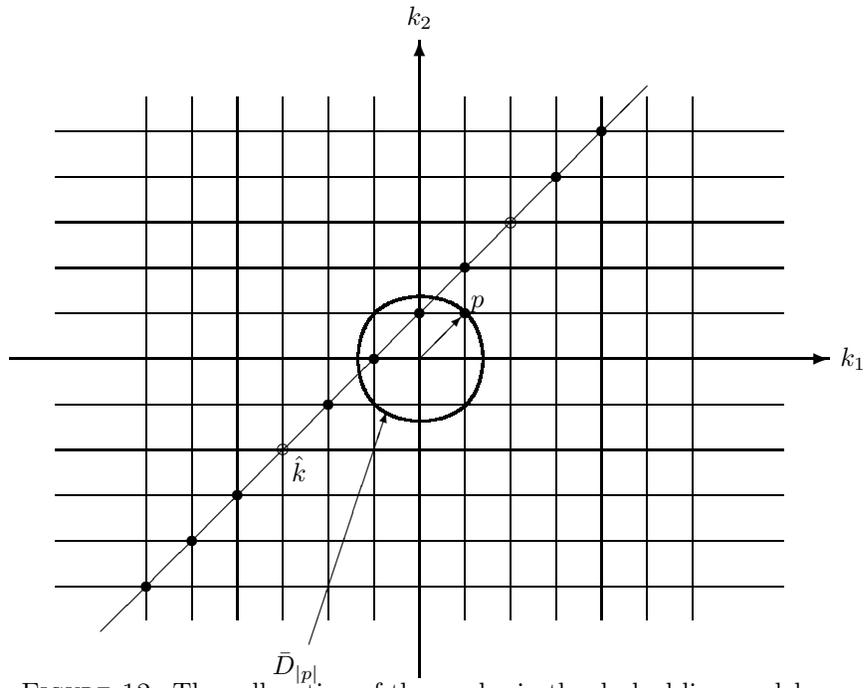
\begin{figure}[ht]
  \begin{center}
    \leavevmode
      \setlength{\unitlength}{2ex}
  \begin{picture}(36,27.8)(-18,-12)
    \thinlines
\multiput(-12,-11.5)(2,0){13}{\line(0,1){23}}
\multiput(-16,-10)(0,2){11}{\line(1,0){32}}
    \thicklines
\put(0,-14){\vector(0,1){28}}
\put(-18,0){\vector(1,0){36}}
\put(0,15){\makebox(0,0){$k_2$}}
\put(18.5,0){\makebox(0,0)[l]{$k_1$}}
\qbezier(-2.75,0)(-2.6375,2.6375)(0,2.75)
\qbezier(0,2.75)(2.6375,2.6375)(2.75,0)
\qbezier(2.75,0)(2.6375,-2.6375)(0,-2.75)
\qbezier(0,-2.75)(-2.6375,-2.6375)(-2.75,0)
    \thinlines
\put(2,2){\circle*{0.5}}
\put(0,0){\vector(1,1){1.85}}
\put(2.275,2.275){$p$}
\put(-12,-10){\circle*{0.5}}
\put(-10,-8){\circle*{0.5}}
\put(-8,-6){\circle*{0.5}}
\put(-6,-4){\circle{0.5}}
\put(-4,-2){\circle*{0.5}}
\put(-2,0){\circle*{0.5}}
\put(0,2){\circle*{0.5}}
\put(2,4){\circle*{0.5}}
\put(4,6){\circle{0.5}}
\put(6,8){\circle*{0.5}}
\put(8,10){\circle*{0.5}}
\put(-14,-12){\line(1,1){24}}
\put(-5.6,-5.4){$\hat{k}$}
\put(-4.4,-13.6){\makebox(0,0)[r]{$\bar{D}_{|p|}$}}
\put(-4.85,-12.55){\vector(1,3){3.4}}
\end{picture}
\end{center}
\caption{The collocation of the modes in the dashed-line model.}
\label{model}
\end{figure}

To simplify our study, we study only the case when $\om_k$ is real, $\forall 
k \in \Z$, i.e. we only study the cosine transform of the vorticity, 
\[
\Om = \sum_{k \in \Z} \om_k \cos (k \cdot X)\ .
\]
To further simplify our study, we will study a concrete 
dashed-line model based upon the line of fixed points (\ref{fixpt}) with the 
mode $p=(1,1)^T$ parametrized by $\Ga$, and the only unstable invariant 
linear subsystem labeled by $\hk = (-3,-2)^T$. The so-called 
{\em{dashed-line model}} is given by \cite{Li01d},
\begin{eqnarray}
\dot{\omega}_n &=& \epsilon_{n-1} A_{n-1} \omega_p
\omega_{n-1} - \epsilon_{n+1} A_{n+1} \omega_p \omega_{n+1} \ , 
\nonumber \\ 
\label{rdlm} \\
\dot{\omega}_p &=& - \sum_{n \in Z} \epsilon_n \epsilon_{n-1}
A_{n-1,n} \omega_{n-1} \omega_n \, , \nonumber
\end{eqnarray}
where
\begin{equation}
  \omega_n = \omega_{\hat{k}+np} \, , \ \  
  A_n = A(p,\hat{k}+np) \, , \ \  
  A_{m,n} = A(\hat{k}+mp,\hat{k}+np) \, ,  
\label{abbn}
\end{equation}
\[
\epsilon_n = \left\{
    \begin{array}{ll}
      1 \, , & \ \ \hbox{if } n \ne 5j \, , \, \forall j \in Z \ ,  \\
\epsilon \, , & \ \ \hbox{if } n = 5j \, , \, 
      \hbox{ for some }j \in Z \, .
    \end{array} \right.
\]
The model is designed to model the hyperbolic structure of the 2D Euler 
equation, connected to 
the line of fixed points (\ref{fixpt}) with $p=(1,1)^T$. Figure \ref{model}
illustrates the collocation of the modes in this model, which has the 
{\em{``dashed-line''}} nature leading to the name of the model. 
\begin{figure}
\caption{The unstable and stable manifolds of the fixed point $\om^*$ 
(\ref{fixpt}): (a) and (b) show two ``painted eggs'' and (c) shows a 
``lip''.}
\label{het}
\end{figure}
When $\e =0$, one can determine the stable and
unstable manifolds of the fixed point $\omega^*$:
\begin{equation}
\omega_p =  \Gamma \, , \, \omega_n=0 \quad (n \in \ZZ) \, .
\label{nfixpt}
\end{equation}
Using the polar coordinates:
\begin{displaymath}
  \omega_1 = r \cos \theta \, , \, 
  \omega_4 = r \sin \theta \, ; \, 
  \omega_2 = \rho \cos \vartheta \, , \, 
  \omega_3 = \rho \sin \vartheta \, ;
\end{displaymath}
we have the following explicit expressions for the stable and
unstable manifolds of the fixed point $\om^*$ and its
negative $-\om^*$ represented through {\em{perverted 
heteroclinic orbits}}:
\begin{eqnarray}
\omega_p &=& \Gamma \ \tanh \tau \, , \nonumber \\
r &=& \sqrt{ \frac{A_2}{A_2-A_1}}\, \ \Gamma \ \mbox{sech}\ \tau \, , \nonumber  \\[1ex]
\theta &=& - \ \frac{A_2}{2\k} \ \mbox{ln}\ \cosh \tau + \theta_0 \, , 
\label{exus}\\[1ex]
\rho &=& \sqrt{\frac{-A_1}{A_2}} \  r \, ,\nonumber  \\[1ex]
\theta + \vartheta &=& \left\{
  \begin{array}{ll}
    - \arcsin \left[ \frac{1}{2} \sqrt{\frac{A_2}{-A_1}}\,  \right] 
        \ , & (\k>0) \ , \\[2ex]
\pi + \arcsin  \left[ \frac{1}{2} \sqrt{\frac{A_2}{-A_1}} \right]
\, , & (\k<0) \, ,
  \end{array} \right.  \nonumber 
\end{eqnarray}
where $\tau = \k \Gamma t 
+ \tau_0$, $(\tau_0, \theta_0)$ are the two parameters
parametrizing the two-dimensional stable (unstable) manifold,
and
\begin{displaymath}
  \k = \sqrt{-A_1 A_2} \cos (\theta + \vartheta) 
   = \pm \sqrt{-A_1 A_2} \sqrt{1+ \frac{A_2}{4A_1}} \ .
\end{displaymath}
The two auxilliary variables $\om_0$ and $\om_5$ have the
expressions:
\begin{eqnarray*}
  \omega_0 &=& \frac{\alpha \beta}{1+ \beta^2} \ \mbox{sech}\ \tau \left\{ 
            \sin [ \beta \ \mbox{ln}\ \cosh \tau +  \theta_0]
            - \frac{1}{\beta} \cos
            [ \beta \ \mbox{ln}\ \cosh \tau + \theta_0 ] \right\} \,
   , \\
  \omega_5 &=& \frac{\alpha \beta}{1+ \beta^2} \ \mbox{sech}\ \tau
  \left\{ \cos [ \beta \ \mbox{ln}\ \cosh \tau + \theta_0 ]
    + \frac{1}{\beta} \sin [ \beta \ \mbox{ln}\ \cosh \tau + \theta_0]
    \right\} \, ,
\end{eqnarray*}
where
\begin{displaymath}
  \alpha = -A_1 \Gamma \k^{-1} \sqrt{\frac{A_2}{A_2-A_1}} \ , \ \ 
  \beta = - \frac{A_2}{2\k} \ .
\end{displaymath}
The explicit expression (\ref{exus}) of the perverted heteroclinic orbits
shows that the unstable manifold $W^u(\om^*)$ 
of the fixed 
point $\om^*$ is the same as the stable 
manifold $W^s(\om^*)$ of the negative $-\om^*$ of 
$\om^*$, and the 
stable manifold of the fixed point $\om^*$ 
is the same as the unstable manifold of the negative $-\om^*$ 
of $\om^*$. Both $W^u(\om^*)$ and $W^s(\om^*)$ have 
the shapes of ``painted eggs'' (Figure \ref{het}). $W^u(\om^*)$ and 
$W^s(\om^*)$ together form the ``lip'' (Figure \ref{het}) which is a 
higher dimensional generalization of the heteroclinic connection on plane.

As the first step toward understanding the degeneracy v.s. nondegeneracy 
of the hyperbolic structure of 2D Euler 
equation, we are interested in the $\e$-homotopy deformation of such 
hyperbolic structure for the dashed-line model. We have calculated 
the Melnikov functions to study the ``breaking'' or ``persistence'' 
of such hyperbolic structure when $\e$ is small. It turns out that both the 
first and the second order Melnikov functions are identically zero 
\cite{Li01d}. 

\section{Conclusion}

Results on chaos in partial differential equations are reported. 
Results on Lax pairs of Euler equations of incompressible inviscid 
fluids are also reported. Preliminary results on dynamical system 
studies of 2D Euler equation are summarized.

\end{document}